\def\pf{{\bf Proof:}}

\def\beg{\begin{equation}}
\def\en{\end{equation}}
\def\beq{\begin{eqnarray}}
\def\enq{\end{eqnarray}}

\def\bbbr{{\rm I\!R}} %reelle Zahlen
\def\bbbn{{\rm I\!N}} %natuerliche Zahlen

\def\bbbc{{\mathchoice {\setbox0=\hbox{$\displaystyle\rm C$}\hbox{\hbox
to0pt{\kern0.4\wd0\vrule height0.9\ht0\hss}\box0}}
{\setbox0=\hbox{$\textstyle\rm C$}\hbox{\hbox
to0pt{\kern0.4\wd0\vrule height0.9\ht0\hss}\box0}}
{\setbox0=\hbox{$\scriptstyle\rm C$}\hbox{\hbox
to0pt{\kern0.4\wd0\vrule height0.9\ht0\hss}\box0}}
{\setbox0=\hbox{$\scriptscriptstyle\rm C$}\hbox{\hbox
to0pt{\kern0.4\wd0\vrule height0.9\ht0\hss}\box0}}}}

\def\bbbg{{\mathchoice {\setbox0=\hbox{$\displaystyle\rm G$}\hbox{\hbox
to0pt{\kern0.4\wd0\vrule height0.9\ht0\hss}\box0}}
{\setbox0=\hbox{$\textstyle\rm G$}\hbox{\hbox
to0pt{\kern0.4\wd0\vrule height0.9\ht0\hss}\box0}}
{\setbox0=\hbox{$\scriptstyle\rm G$}\hbox{\hbox
to0pt{\kern0.4\wd0\vrule height0.9\ht0\hss}\box0}}
{\setbox0=\hbox{$\scriptscriptstyle\rm G$}\hbox{\hbox
to0pt{\kern0.4\wd0\vrule height0.9\ht0\hss}\box0}}}}

\def\bbbq{{\mathchoice {\setbox0=\hbox{$\displaystyle\rm Q$}\hbox{\raise
0.15\ht0\hbox to0pt{\kern0.4\wd0\vrule height0.8\ht0\hss}\box0}}
{\setbox0=\hbox{$\textstyle\rm Q$}\hbox{\raise
0.15\ht0\hbox to0pt{\kern0.4\wd0\vrule height0.8\ht0\hss}\box0}}
{\setbox0=\hbox{$\scriptstyle\rm Q$}\hbox{\raise
0.15\ht0\hbox to0pt{\kern0.4\wd0\vrule height0.7\ht0\hss}\box0}}
{\setbox0=\hbox{$\scriptscriptstyle\rm Q$}\hbox{\raise
0.15\ht0\hbox to0pt{\kern0.4\wd0\vrule height0.7\ht0\hss}\box0}}}}
\def\bbbt{{\mathchoice {\setbox0=\hbox{$\displaystyle\rm
T$}\hbox{\hbox to0pt{\kern0.3\wd0\vrule height0.9\ht0\hss}\box0}}
{\setbox0=\hbox{$\textstyle\rm T$}\hbox{\hbox
to0pt{\kern0.3\wd0\vrule height0.9\ht0\hss}\box0}}
{\setbox0=\hbox{$\scriptstyle\rm T$}\hbox{\hbox
to0pt{\kern0.3\wd0\vrule height0.9\ht0\hss}\box0}}
{\setbox0=\hbox{$\scriptscriptstyle\rm T$}\hbox{\hbox
to0pt{\kern0.3\wd0\vrule height0.9\ht0\hss}\box0}}}}
\def\bbbs{{\mathchoice
{\setbox0=\hbox{$\displaystyle     \rm S$}\hbox{\raise0.5\ht0\hbox
to0pt{\kern0.35\wd0\vrule height0.45\ht0\hss}\hbox
to0pt{\kern0.55\wd0\vrule height0.5\ht0\hss}\box0}}
{\setbox0=\hbox{$\textstyle        \rm S$}\hbox{\raise0.5\ht0\hbox
to0pt{\kern0.35\wd0\vrule height0.45\ht0\hss}\hbox
to0pt{\kern0.55\wd0\vrule height0.5\ht0\hss}\box0}}
{\setbox0=\hbox{$\scriptstyle      \rm S$}\hbox{\raise0.5\ht0\hbox
to0pt{\kern0.35\wd0\vrule height0.45\ht0\hss}\raise0.05\ht0\hbox
to0pt{\kern0.5\wd0\vrule height0.45\ht0\hss}\box0}}
{\setbox0=\hbox{$\scriptscriptstyle\rm S$}\hbox{\raise0.5\ht0\hbox
to0pt{\kern0.4\wd0\vrule height0.45\ht0\hss}\raise0.05\ht0\hbox
to0pt{\kern0.55\wd0\vrule height0.45\ht0\hss}\box0}}}}
\def\bbbz{{\mathchoice {\hbox{$\sf\textstyle Z\kern-0.4em Z$}}
{\hbox{$\sf\textstyle Z\kern-0.4em Z$}}
{\hbox{$\sf\scriptstyle Z\kern-0.3em Z$}}
{\hbox{$\sf\scriptscriptstyle Z\kern-0.2em Z$}}}}

\def\gaaa{\ifmmode
              {{\mbox{\deu a}}}%
          \else${{\mbox{\deu a\ }}}$%
           \fi}
\def\gbbb{\ifmmode
              {{\mbox{\deu b}}}%
          \else${{\mbox{\deu b\ }}}$%
           \fi}
\def\ccc{\ifmmode
        {\bbbc}%
          \else${\bbbc\ }$%
          \fi}
\def\fff{\ifmmode
        {{\mbox{\bf F}}}%
          \else${{\mbox{\bf F\ }}}$%
          \fi}
\def\fxfx{\ifmmode
              {{\bbbf}_p[x] }%
          \else${{\bbbf}_p[x]\ }$%
          \fi}
\def\fffp{\ifmmode
              {{\mbox{\bf F}}_p }%
          \else${{{\mbox{\bf F}}}_p\ }$%
          \fi}
\def\fffq{\ifmmode
              {{\bbbf}_q }%
          \else${{\bbbf}_q\ }$%
          \fi}
\def\gmmm{\ifmmode
              {{\mbox{\deu m}}}%
          \else${{\mbox{\deu m\ }}}$%
           \fi}
\def\nnn{\ifmmode
        {\bbbn}%
          \else${\bbbn\ }$%
          \fi}

\def\nnno{\ifmmode
        {{\bbbn}_0}%
          \else${{\bbbn}_0\ }$%
          \fi}

\def\ooo{\ifmmode
              {{\mbox{\deu o}}}%
          \else${{\mbox{\deu o\ }}}$%
           \fi}
\def\gppp{\ifmmode
              {{\mbox{\deu p}}}%
          \else${{\mbox{\deu p\ }}}$%
           \fi}

\def\ppp{\ifmmode
              {{\mbox{\deu p}}}%
          \else${{\mbox{\deu p\ }}}$%
           \fi}
\def\cppp{\ifmmode
              {\cal P}%
          \else${\cal P\ }$%
           \fi}

\def\qqq{\ifmmode
              {\bbbq}%
          \else${\bbbq\ }$%
           \fi}
\def\zzz{\ifmmode
              {\bbbz}%
          \else${\bbbz\ }$%
          \fi}
\def\grrr{\ifmmode
              {{\mbox{\deu r}}}%
          \else${{\mbox{\deu r\ }}}$%
           \fi}
\def\rrr{\ifmmode
              {\bbbr}%
          \else${\bbbr\ }$%
          \fi}
\def\st{\ifmmode
              {***}%
          \else${***}$%
          \fi}

\newcommand{\deu}{\bf}

\newcommand{\text}{\mbox}

\documentstyle[12pt]{article}
\begin{document}
\title{ Exceptional linear systems on curves on Enriques surfaces}
\author{Severinas Zube}
\date{}
\maketitle

\begin{abstract}
The main purpose in this paper is to study the gonality, the Clifford index and
the Clifford dimension on linearly equivalent smooth curves on Enriques
surfaces. The method is similar to techniques of M.Green $\&$ R.Lazarsfeld and
G.Pareschi.
\end{abstract}

\subsection{ Introduction}
In recent years several authors have been led to study the following question :
to what extent do linearly equivalent smooth curves on a surface $S$ carry
"equally exceptional" linear series? Green and Lazarsfeld investigated the case
when $S$ is a K3 surface and  proved in [GL] that smooth linearly equivalent
curves have the same Clifford index . Clifford index is a natural numerical
invariant measuring "exceptional" linear systems on curves . On the other hand
Donagi's example shows that smooth linearly equivalent curves on K3 surfaces
can have diferent gonality. The question in the case of Del Pezzo surfaces of
degree $\geq 2$ has been studied  by Pareschi . It follows from the result in
[P] that gonality and Clifford index is the same  for all smooth linearly
equivalent curves, with one exception for gonality , involving curves of genus
3 (see [P] for details  ).

The purpose of this note is to study the same question in the case of Enriques
surfaces. It turns out that the gonality and the Clifford index is not constant
for smooth linearly equivalent curves (see section 5.11 ) and in general we
obtain  the following estimate for the jump of the gonality and the Clifford
index in  a linear system:

\newtheorem{ttt}{Theorem}
%[subsection]
\begin{ttt}
Let C be a smooth irredecible curve of genus g on an Enriques surface S and let
gon($|C|$) = min \{gon($C'$) $|$ $C'$ is smooth curve in $|C|$ \}, cliff($|C|$)
= min \{cliff($C'$) $|$ $C'$ is smooth curve in $|C'|$ \} be the minimal
gonality , Clifford index for $|C|$ . Then for all smooth curves $C'$ in $|C|$
gonality $C'$ is $ \leq gon(|C|) +2$. And if cliff($C'$) = $gon(C') -2$ then
Clifford index of $C'$ is $\leq cliff(|C|) +2 $. Moreover, if $gon(|C|) \leq
\frac{g-1}{2}$ or $g \geq9$ , then there is a line bundle L on S such that
$cliff(C') \leq cliff(L|_{C'} ) +2$ .( see proposition~\ref{PPP} for more
precise result ).
\end{ttt}

\paragraph{Remark}: \\
- It seems  that Cliford index always can be computed by gonality, i.e.
$cliff(C)$  $= gon(C) -2$ ,  for all curves contained in Enriques surfaces.
(see Theorem 2 and conjecture bellow ). \\
-In the section 5.11 we give examples which show that the gonality for the
smooth linearly equivalent cuves is not constant. \\

 Another very useful invariant  describing "exceptional" curves is the Clifford
dimension introduced by Eisinbud , Lange , Martens and Schreyer in [ELMS]. If
the curve has the Clifford dimension greater than 1 then it is "special" (see
[ELMS] for description of such curves ). For the curves contained in $S$ we
have the following result:

%\addtocounter{subsection}{1}
\newtheorem{cliff}[ttt]{Theorem}
%[subsection]
\begin{cliff}
 Let C be a smooth irredecible curve on an Enriques surface S then its Clifford
 dimension is 1 or greater than 9. ( In fact we prove a somewhat more precise
result - see Proposition~\ref{CCC} ).
\end{cliff}

We conjecture that every   smooth irredecible curve $C$ on Enriques surface S
has Clifford dimension 1.

The paper is organized as follows: In the Section 1.? we introduce notation and
collect some preliminary results. In Section 2.? we investigate the bundle
$E(C,A)$ - introduced by Lazarsfeld - and obtain the jumping estimate for
gonality and prove Theorem 1. Section 3.? is devoted to the proof that a smooth
 plane curve  of the degree $\geq5$ cannot be contained in an Enriques surface
. In section 4.? we prove Theorem 2. In last section we explain how to obtain ,
using results of section 2.? , explicit examples of the smooth linearly
equivalent curves of different gonality .

This paper was written during my stay at SFB 170 in  G$\ddot{o}$ttingen. I am
grateful SFB 170 for support and good working condition.

\paragraph{Notation}: \\
-We will work over the complex field. \\
-We denote by the same capital letter both a divisor and the line bundle
associate to the divisor and we hope that the meaning  will be clear from
context.

 %\paragraph{Preliminaries} : \\
\subsection*{Preliminaries}
\addtocounter{section}{1}
%\addtocounter{subsection}{-3}
\subsection{E(C,A)}
Let $S$ be a regular surface , $C$ a curve contained in $S$, and $A$ a base
point free line bundle on $C$. One can associate to the triple $(S,C,A)$ a
certain vector bundle on $S$ in a canonical way. We refer the reader to [P],
[GL], [T] for  details.

Let us denote by $F(C,A)$ the vector bundle defined by the sequence:
$$0 \to F(C,A) \to H^{0}(A) \otimes O_S \stackrel{ev}{\to}A \to 0, $$   where
$A$ is viewed as a sheaf on $S$ and $ev$ is the evaluation map.

Dualizing the above sequence we get (see [P])
 \begin{equation} %{subsection}
0 \to H^{0}(A)^{*} \otimes O_S \to E(C,A) \to O_{C}(C) \otimes  A^ * \to 0
\label{E}
\end{equation}
here $E(C,A)=F(C,A)^ * = {\cal H}{\sl om}(F(C,A),O_S )$ , $A^* = {\cal H}{\sl
om}(A,O_C )$  . To simplify the notation  we will omit  $(C,A)$, if it is clear
from context to which pair $(C,A)$ the bundle $E(C,A)$ is associated. We have
\begin{equation}
c_1 (E)=C; {~~~}
c_2 (E) = deg (A);
\end{equation}
\begin{equation}
  rk(E)=h^0 (A);{~~~~}
h^i (F(C,A)) =h^{2-i}(E \otimes K_S) = 0{~~~} for {~~~~} i=1,2;
\end{equation}
where $K_S$ is the canonical bundle of $S$. If $h^0 (O_{C}(C) \otimes  A^ * ) >
0$, then  $E(C,A)$  is generated by its sections away from a finite set
coinciding with the fixed divisor of $O_{C}(C) \otimes  A^ *   $

Let $s^ \perp $ be a subspase of $H^0 (A)^* $ orthogonal to $ s \in H^0 (A) $
with respect to  the  natural pairing
$$ H^1 (K_C -A) \otimes H^0 (A) = H^0 (A)^* \otimes H^0 (A) \to H^1 (K_C )=
\ccc.$$
Then we get another description of the bundle $E(C,A)$:
\begin{equation}
0 \to s^ \perp \otimes K_S \to E(C,A) \otimes K_S \to J_ \xi (C+K_S ) \to 0,
\label{ej}
\end{equation}
 where $ \xi = (s)_0$ is the zero scheme of the section $s \in H^0 (A) $ (see
[T]).

%\addtocounter{section}{1}
\subsection{Enriques surfaces}

A smooth irreducible surface $S$, such that $h^1 (O_S )=h^2 (O_S )=0 $ and
$2K_S \sim O_S $, is called a Enriques surface.   If $C$ is a smooth curve in
$S$ , then by adjunction formula one has :
\beg
g(C) = \frac{C^2 }{2} +1 . \label{g}
\en
Recall that a divisor $D$ on a smooth surface $X$ is said to be $nef$ if $DC
\geq 0$ for every curve $C$ on $X$ . The following properties will be  used
throughout, sometimes without explicit mention:

(A)([C,D] Proposition 3.1.6) If $D$ is a $nef$   divisor , then $ \mid D \mid $
has no fixed components , unless $  D \sim  2E+R $, where $ \mid 2E \mid $ is a
genus 1 pencil and $R$ a smooth rational curve with $RE=1$ . \label{111}

(B)([C,D] Corollary 3.1.3) If $D$ is a $nef$ divisor  and $D^2 > 0 $ , then
$H^1 (O_S (-D)) = 0$ and $\chi(O(D)) -1 =dim|D| = \frac{D^2 }{2}$. \label{222}

(C)([C,D] Proposition 3.1.4) If $ \mid D \mid $ has no fixed components,   then
one of the following holds: \nonumber \\
(i) $D^2 > 0$ and there exist an irreducible curve $C$ in $ \mid D \mid $.
\nonumber \\
(ii) $D^2 =0$ and there exist a genus 1 pencil $ \mid P \mid $ such that $D
\sim kP$ for some $k \geq 1 $. \label{333}

(D)([C,D]Chapter 4, appendix , corollary 1. and corollary 2.) If $D^2 \geq 6$
and $D$ is $nef$ then $D$ is ample , $2D$ is generated by its global sections,
$3D$ is very ample.

Let $K^+ = \{D^2 > 0 | D {~}is{~} an {~}effective{~} divisor. \}$. $ K^+ $ is
called the positive cone.     $ \overline{K^+ } = \{D^2 \geq 0 | D{~} is
{~}an{~} effective{~} divisor. \}$    is  the closure of the positive cone.

The Enriques surface $S$ is called unnodal if there are no smooth (-2) curves
contained in $S$.

%\addtocounter{section}{1}
\subsection{Steiner construction}

Let $L,M$ be the effective divisors on $S$. Parametrize the two pencils $ \mid
L( \lambda )  \mid  \\ \subset \mid L \mid $ and $ \mid M( \lambda ) \mid
\subset \mid M \mid $ by $ \lambda \in P^1 $, choosing the parametrization so
that $L( \lambda ) $ and $ M( \lambda ) $ have no common components for every $
\lambda \in P^1 $.

Then the curve :
\beg
C= { \cup }_{ \lambda} (L( \lambda ) \cap M( \lambda ))
\en
is cleary irreducible if the general curves in $ \mid L( \lambda ) \mid  $ and
$ \mid M( \lambda ) \mid   $  are irreducible. Denote by $ \xi $ (resp. $\eta
$) base locus of $ \mid L( \lambda )  \mid $ (resp.  $ \mid M( \lambda ) \mid
$).  Then $ \mid L( \lambda ) \mid  =   \mid J_{ \xi } (L) \mid , \mid M(
\lambda ) \mid =  \mid J_{ \xi} (M) \mid $. If $ \xi \cup \eta = P_1 + ...
+P_{L^2+M^2}$ consist of different points , then $C$ contain $ \xi \cup \eta $.
We have $C \sim L+M $ . Indeed , using our data we  can obtain the  sequence:
\beg
0 \to \ccc^2 \otimes O_S \stackrel{ s}{ \to} M \oplus L \to { \cal L} \to 0
\en
 where $s(( \mu, \nu ) \otimes f ) = s^L _{ \mu  / \nu }(f) \oplus s^M _{ \mu
/  \nu }(f) $ and the map $s^L _{ \mu  /  \nu } ( \bullet ) {~}$ (resp. $s^M _{
\mu  /  \nu } ( \bullet ) $ ) is  multiplication  by $ L( \mu  /  \nu ) \in|L|
$ (resp.$ M( \mu  /  \nu) \in|M|$), hence is defined by the sequences:
\beq
0 \to O \to L \to O_{L( \mu  /  \nu )} \to 0 & (s^L _{ \mu / \nu } )\nonumber
\\
(resp.{~} 0 \to O \to M \to O_{M( \mu  /  \nu )} \to 0 & (s^M _{ \mu / \nu })
{~}). \nonumber
\enq
We see that the sheaf  $ \cal L $  has support  $C$ and $c_1 (M \oplus L) \sim
M+L \sim supp { \cal L} = C$.

In the situation described above we will say that $C$ is obtained by the
Steiner constrution using $ \mid L( \lambda) \mid $ and $ \mid M( \lambda) \mid
$ .

Note , that if a curve $C' \in \mid C \mid $ contains the base locus of $ \mid
L ( \lambda ) \mid $ , then $C'$ can be obtained by the Steiner construction
using $ \mid L ( \lambda ) \mid $ and any another pencil $ \mid M' ( \lambda )
\mid   \subset \mid M \mid $ such that $  L ( \lambda ) $  and  $ M' ( \lambda
)  $
have no fixed components for every $ \lambda \in P^1 $.

%\addtocounter{section}{1}
\subsection{Gonality , Clifford index}

Let $C$ be a curve of genus $g$ and and $A$ be a line bundle on $C$ . The
Clifford index of $A$ is the integer
\[cliff(A) = deg(A) -2r(A) ,\]
where $ r(A) = h^0 (A) -1$ is the projective dimension of $|A|$. \\
The Clifford index of $C $ is
\[cliff(C) = min \{ cliff(A) {~}|{~}r(A) \geq 1, {~} deg(A) \leq g-1 \}.\]
Note that this last definition makes sense only for curves of genus $g \geq 4$.

We say that a line bundle $A$ contributes to the Clifford index  if it
satisfies the inequalities in the above definition and that $A$  computes the
Clifford  index  if $A$ contributes to the Clifford index and $cliff(C) =
cliff(A)$.
Finaly, we define the Clifford dimension of $C$ as
\[ cliffdim(C) = min \{r(A) , A {~}computes {~}the {~}Clifford {~}index{~}of
{~}C \}. \]

The gonality of $C$, denoted by $gon(C)$, is the minimal degree of a pensil on
$C$. Such a definition is nontrivial only  for curves of genus $g \geq 3$. Also
we say $A$ computes to the gonality if $r(A) = 1$ and $gon(C) = deg(A) $. By
Brill-Noether theory one has that
\[cliff(C) \leq gon(C) -2 \leq \left[ \frac{g-1}{2} \right] \]
and for a general curve both inequalities are equalities.

\addtocounter{section}{1}
\subsection{Gonality of curves on Enriques surfaces}
We will denote (see [C,D] p. 178 ) by
\[ \Phi:K^+  \to \bf{Z} _{ \geq 0} \]
the function defined by
\[ \Phi(C) = inf \{ CE, |2E| {~} is {~}a{~}genus{~} 1 {~}pencil \} . \]

\newtheorem{lll}{Lemma }[subsection]
\begin{lll}
Let $C$ be a smooth curve of genus $g$ on an Enriques surface $S$. Then: \\
(i) $  \Phi(C) \leq  \left[ \sqrt{2g-2} \right]$ , where $[l]$ means the
integer part of $l$. \\
(ii) if $2 \Phi (C) \leq g-1$  and $  \Phi (C) = CE$ ($ |2E| $ is a genus 1
pencil) , then $cliff(C) \leq  cliff(2E|_C ) \leq 2 \Phi(C) - 2$ and $2E|_C$
contributes to the Clifford index. ( by (i) condition $2 \Phi(C) \leq g-1$ is
always satisfied if $g \geq 9$). \label{lll}
\end{lll}
\pf \ \
(i)  \ \ is a statment about  the Enriques lattice $H^2 (S, \bf{Z})$ and has
been proved in [C,D] (Corollary 2.7.1). \\
(ii)  Assume $ \Phi (C) = CE $ , where $|2E|$ is a genus 1 pencil. Then we have
the sequence:
\[ 0 \to O(2E-C) \to O(2E) \to O(2E)|_C \to0 . \]
 The divisor $2E-C$ cannot be an effective, because $(2E-C)E < 0$ . Therefore
$h^0 (2E|_C ) \geq 0 $ and $cliff(C) \leq  cliff(2E|_C ) \leq 2 \Phi (C) -2 $,
since by our condition $2E|_C$ contributes to the Clifford index.$ \Box$

\newtheorem{PPP}{Proposition}[subsection]

\begin{PPP}
 Let $C$  be a smooth curve on an Enriques surface $S$, and
let $A$ be a pencil on $C$ such that $deg(A)=gon(C) \leq  \frac{g-1}{2}$. Then

(i) There exist an exact sequence:
\beg
0 \to M \to E(C,A) \to L \to 0 , \label{mle}
\en
such that M, M-L are line bundles from  the positive cone $K^+$   .  Moreover
the linear systems $\mid L \mid , \mid M \mid $ have no fixed components
and              $M^2 > LM = deg(A) > L^2 \geq 0$  (if $deg(A) = \frac{g-1}{2}
$  then  $M-L \in \overline{K^+}$ and $M^2 \geq LM = deg(A) \geq L^2 \geq 0$).

(ii) The curve C can be obtained by the Steiner construction for some

 {~~~~~}$|M( \lambda)| = P^1 \subset \mid M \mid {~~}and {~~}|L( \lambda)| =
P^1 \subset \mid L \mid $.

(iii) If S is an unnodal Enriques surface , then

$ {~~} $     (a) The exact sequence ( \ref{mle}) splits.

$ {~~} $     (b) If $gon( \mid C \mid)= gon(C)$ then $L \sim 2E_1  {~} or {~}
E_1 + E_2 $, where $ \mid 2E_1 \mid $,  \\ $\mid 2E_2 \mid     $ are two genus
1 pencils  on S. \label{PPP}
\end{PPP}
\pf  \ \ \   1 Case. \ \ \   $d< \frac{g-1}{2}$.

(i). Assume that $d < \frac{g-1}{2} $ , then the vector bundle $E(C,A)$ is not
stable in the Bogomolov sense. Indeed $4c_2 (A) = 4deg(A) \leq c_1 (E(C,A) =
C^2=2g-2$.  By Bogomolov`s theorem [B] we have an exact sequence:
\beg
0 \to M \to E(C,A) \to J_{\xi}(L) \to 0 , \label{new}
\en
 that such $ (M-L)^2 > 0 $ and $M-L$ is an effective divisor.  $E(C,A)$ is a
vector bundle  genereted by its global sections away from a finet set.
Therefore $J_ \xi (L)$ is also generated by its global sections away from a
finet set and we obtain that $L^2 \geq 0$ and linear system $ \mid L \mid $ has
not fixed components.

Claim: \ \  $l( \xi)= length( \xi) = 0$.

If $L^2= 0$, then by proposition 3.14 [CD] $L \sim kP$ for some $k \geq 1$.
By assumtion, $ deg(A) = kPM+l( \xi) = gon(C)$, but $deg(P \mid _C) =PM \geq
deg(A)$, hence we have $l( \xi) = 0, {~}k=1$ .

Assume $L^2 > 0$ and $l( \xi) > 0 $, then from exact sequence (~\ref{new} ) we
obtain:
\begin{eqnarray}
h^0 (E(-L)) &=& h^0 (M-L) , \label{0} \\
h^1 (E(-L)) &=& l( \xi) + h^1 (M-L) , \label{H22} \\
h^2 (E(-L)) &=& 0.
\end{eqnarray}

On the other hand by 1.3 (B) $H^1 (-L) =0$, therefore using (~\ref{E} ) we
have:
\beq
h^0 (E(-L)) &=& h^0 (C,M \mid _C -A) , \label{13} \\
h^1 (E(-L)) &=&   h^1 (C,M \mid _C -A) -2h^2  (-L) , \label{H11} \\
h^2 (E(-L)) &=& 0.
\enq

By Riemann-Roch on the curve C

\beq
 h^0 (C,M \mid _C  -A)  =  MC-A-g(C)+1+h^1 (M \mid _C -A)   &=&
(by{~}~\ref{H11} ,~\ref{H22}{~} and {~} duality){~~}  \nonumber  \\
   =MC - A -( \frac{M^2 + L^2}{2} + ML +1)+1 + \nonumber \\
+ l( \xi) +2h^0 (L \otimes K_S )+ h^1 (M-L) &=& (by {~~} R-R {~~} and {~~} (
\ref{new} ) )  \nonumber \\
= M^2   -( \frac{M^2+L^2}{2} + ML +1)+ L^2 +2+ h^1 (M-L)  &= &  \nonumber \\
 = \frac{M^2+L^2}{2} - ML +2+ h^1 (M-L). \label{l}
\enq

We have $h^2 (M-L) = 0$, hence (~\ref{0}), (~\ref{13} ) and (~\ref{l} ) gives
us
\beg
\chi (M-L)= \frac{M^2+L^2}{2} - ML +2 \label{17}
\en

 But by Riemman-Roch theorem $\chi(M-L) = \frac{M^2+L^2}{2} -ML+1$. This
contradicts (~\ref{17} ), hence $l( \xi)=0 $ .

Now we see that $L$ is nontrivial , because $c_2 (E) = deg(A) = LM$. This
implies
$(M-L)L > 0$ and $(M-L)M > 0$ , since $ M-L,M \in K^+;L \in \overline{ K^+}$.

(ii). By ( \ref{ej}) we have:
\beg
0 \to O_S \to E(C,A) = M \oplus  L  \to J_ \xi (C ) \to 0, \label{w}
\en
where $ \xi = (s)_0 $ is the zero scheme of the section $s \in H^0 (A)$ . Since
$H^1 (M) = 0$ we see that  the sequence ( \ref{mle} ) is exact on a section
level i.e.
\[0 \to H^0 (M) \to H^0 (E(C,A)) \to H^0 (L) \to 0 . \]
Hence we can write $H^0 (A) \ni s=s_1 \oplus s_2 $ , where $s_1 \in H^0 (L),
s_2 \in H^0 (M)$. When $s$  runs through  $H^0 (A) $ ,  $s_1 $ and $s_2$ runs
through $|L( \lambda)| = P^1 \subset \mid L \mid$ and $|M( \lambda)| = P_1
\subset \mid M \mid $ respectively. The zero set $ \xi $ consists of  points
where  both sections $s_1 ,s_2$ are zero. This shows  that $C$ is obtained by
the Steiner construction using $|M( \lambda)| ,  |L( \lambda)|$ . Now we can
see , that the linear system $ \mid M \mid  $ has not fixed components, because
a fixed part of $ \mid M \mid$ should be also the fixed part of $ \mid C \mid$
, but the curve $C$ is smooth.

(iii), \ \  (a).
By 1.3 (B)  we see that the extention group $Ext^1 (L,M) = H^1 (M-L)=0$,
therefore the sequence ( \ref{mle} ) splits.

(iii), \ \   (b) . By reducibility lemma 3.2.2 [CD], $L$ is lineary equivalent
to a  sum of genus 1 curves. Assume $L \sim L_1+L_2$ , where $L_1, L_2$ are
effective nontrivial divisors on $S$ such that $L_2$ has no fixed components.
Then a curve $C'$ $ \in \mid C \mid $ which is obtained by the Steiner
construction using some $|L_2 ( \lambda)| \subset   \mid L_2 \mid $ and
$|M+L_1 ( \lambda)| \subset \mid M+L_1 \mid$ has gonality $(M+L_1 )L_2  $ . But
this number is smaller  than $ML= (M' + L)(L_1 + L_2 )= ML_2 + M'L_1 + LL_1$ ,
where $M'= M-L$. This contradicts our assumtion about minimality of $gon(C)$ .
Therefore , there are no such spliting $L$ into two bundles $L_1 ,L_2 $ and we
get (iii) , (b).

2 Case. \ \ \  $ d = \frac{g-1}{2}$.

Assume $ d = \frac{g-1}{2}$, then by the Riemann-Roch theorem :
\beg
 \chi (E(C,A),E(C,A)) = 4+c_1 ^2(E) - 4c_2(E) = 4. \nonumber
\en
If $E(C,A)$ is $H$-stable in the sense of Mamford-Takemoto ($H$ is an ample
divisor on $S$) , then it is simple and $Ext^2 (E,E)^*  = Hom(E,E \otimes K_S )
= \ccc $ or 0. Indeed , any nontrivial homomorphism between $E$ and $E \otimes
K_S $ should be an isomorphism , because both bundles are stable and have the
same determinant. This contadiction shows that $E$ is not $H$-stable for any
ample divisor $H$. By 1.3 (D) $C$ is a ample and hence $E$ is not $C$-stable ,
therefore we have the sequence:
 \[ 0 \to M \to E(C,A) \to J_{ \xi} (L) \to 0 , \]
 such that $(M-L)C = M^2 -L^2  \geq 0 $ . Sence $4(ML+l( \xi ))= 4c_2 (E) =c_1
(E) = M^2 + L^2 + 2ML $ we have $ (M-L)^2 =4l( \xi ) \geq  0$ \ i.e. \  $ M-L
\in \overline{K^+ }$. And now we can argue as in case 1. above. $ \Box$

{\bf Proof of Theorem 1 } : If $gon(|C|) \leq \frac{g-1}{2}$, then by
proposition~\ref{PPP} we are done. On the other hand by Brill-Neother theory we
have :
\[gon(C) \leq \left[ \frac{g-2}{2} \right] +2 . \Box \]

\addtocounter{section}{1}
\subsection{Curves of the Clifford dimension 2}

It is well known that the curve of Clifford dimension 2 is a smooth plane curve
of degree $d \geq 5$ and a line bundle $A$ computing the Clifford index is
unique. In this case there is a 1-dimensional family of pencils of degree $d-1$
computing gonality , all obtained by projecting from a point of the curve.

\newtheorem{ggg}{Proposition}[subsection]
\begin{ggg}
An Enriques surface  does not contain any smooth plane curve of degree $d \geq5
$ .
\end{ggg}

\pf \ \ Let $C$ be a  curve  of degree $ d \geq 6 $ . Recall that for smooth
plane curves of degree $d$ by adjunction formula we have $g(C)
=\frac{d(d-3)}{2} +1$.

 By lemma~\ref{lll} (i)
\[ \Phi(C) \leq \left[ \sqrt{2g(C)-2} \right] = \left[ \sqrt{d(d-3)} \right]
\leq d-2 .\]
On the other hand we have
\[  \Phi(C) \leq d-2 \leq \frac{d(d-3)}{4} =  \frac{g(C) -1}{2}. \]
By lemma~\ref{lll} (ii) we obtain
\[ cliff(C) =d-4 \leq cliff(2E|_C) \leq 2 \Phi(C) -2 \leq d-4 \]
where $\Phi(C) =CE $ and $|2E|$ is a genus 1 pencil . Moreover , $(2E)|_C $
contributes to the Clifford index . This is a contradiction , since we obtain a
pencil which computes the Clifford index.

If $deg(C)= 5$, then from the exact sequence (~\ref{ej} ) we have :
\[h^0 (J_{ \xi} (C+K)) = h^0 (K_C -A) = h^1 (A) = g-d+r(A) = 3, \]
where $\xi = (s)_0 $ is zero set of the section $s \in H^0 (A) $ . Therefore we
can consider $\xi $ as a divisor on $C$ which is linear equivalent to $A$. A
curve from the linear system $D \in |J_A (C+K)|$ cut out   a divisor $A+A_D$ on
our curve $C$ . The divisor $A+A_D$  is linearly equivalent to the canonical
divisor $K_C $, hence $deg(A_D ) = 5$ . $D $ runs through $ |J_A (C+K)|$ and
cuts out  $A_D$ on the curve $C$,  therefore $h^0 (C,A_D) \geq h^0 (C,A) =3$.
By the Cliford theorem $deg(A_D ) \geq 2(h^0 (C,A_D )-1)$ , therefore $h^0
(C,A_D ) = 3$.  The line bundle $A$ is  unique linear systerm of degree 5 and
projective dimension 2 , hence we have $A \sim A_D$ . Denote by $C'$ a smooth
curve in the linear system $|J_A (C+K)|$ which cuts out the divisor $A+A'$ on
our curve $C$ . Then we can consider the divisor $A'$ as a divisor on $C'$. A
curve  $D' \in |J_A (C+K)|$ cuts out a divisor $A+A_{D'}$ on the curve $C'$
($A_{D'} \in|A'|$). In the same way as above we get  $de!
g(A_{D'}) =5 , h^0 (C',A_{D'}) =3$
% Indeed, $|J_{A'} (C)= 3 $

\addtocounter{section}{1}
\subsection{ Curves of the Clifford dimension $r \geq 3$}

The curves of Clifford dimension $r \geq3$ are extremly rare . We refer the
reader to the
preprint [ELMS] for the details. The main result in [ELMS] is: \\
If $C$ has Clifford dimension $r \geq 3$ then one of the following holds:

1. (i) $C$ has genus $g=4r-2$ and Clifford index $2r-3 =[ \frac{g-1}{2}] -1$ .
\\
  \  \ \ \ (ii) $C$ has a unique line bundle $L$ computing the Clifford index ,
$L^2$ is  the canonical bundle on $C$.

2. If the curve does not satisfy  1. , then $ r \geq10$ and the  degree $d$ of
the  line bundle computing the Clifford index is $ \geq 6r-6$ and its genus  is
$ \geq 8r-7$. Eisinbud , Lange , Martens and Schreyer  conjecture  that such a
curve does not exists (see [ELMS]).

\newtheorem{CCC}{Proposition}[subsection]
\begin{CCC}
Let $C$ be a smooth curve of  Clifford dimension $r \geq 3$ and genus $g$
contained in an  Enriques surface $S$. Then $C$ is the curve described by the
case 2 above. Moreover $g \geq 2(r-1)^2  +1 \geq 163$ and $2 \left[ \sqrt{2g-2}
\right] -2 +2r \geq d $, where $d$ is the degree of a line bundle computing
Clifford index. \label{CCC}
\end{CCC}
\pf \ \ Consider the curve which satisfies 1.  above . By lemma~\ref{lll} (ii)
we have
\[cliff(C) \leq 2 \Phi(C) -2 \leq \frac{g-1}{2} -2 .\]
This  contradicts the condition $cliff(C) = [ \frac{g-1}{2}] -1 $.

For a curve which satisfies 2  lemma~\ref{lll} (i) and (ii) implies that
\[ d-2r \leq 2 \Phi(C) -2 \leq 2 \left[ \sqrt{2g-2} \right] -2 .{~~~~~~~}  (*)
\]
But $d \geq 6r-6$ , hence
\[ 4r-6 \leq 2[ \sqrt{2g-2}] -2 \]
and we get
\[ g \geq 2(r-1)^2 +1 \geq 163 . \]
 Also from (*) we obtain
\[d \leq2r-2  + 2 \left[ \sqrt{2g-2} \right] . \Box \]

\addtocounter{section}{1}
\subsection{ Examples}

In this section we will apply proposition~\ref{PPP} to obtain examples of
smooth linearly equivalent curves of different gonality. In explicit examples
we try to explain the reasons why the gonality is not constant .

Assume $S$ is an unnodal Enriques surface. Consider $L = E_1 +E_2 , M = 2L $,
where $|2E_1 | , |2E_2 |$ are genus 1 pencils such that $E_1 E_2 = 1$ , then
$C \sim M+L$ is very ample by corollary 2 ( appendix , after chapter 4 ) [CD].

If $C$ is a smooth curve in the linear system $|M+L|$  which is obtained by
Steiner construction using $P^1 =|L|$ and $|M( \lambda)| \subset |M|$ , then by
proposition~\ref{PPP}  we see that $gon(C) \leq  ML =4$  and actually $gon(C) =
4$. Indeed, if $gon(C) \leq 3$ then by proposition~\ref{PPP} we obtain  $M'$
and $L'$ such that $C \sim L'+M' , gon(C) =M'L' > L'^2 \geq 0 $ . Therefore
$L'^2 \leq2$ and  by (iii)(b) we see that $L \sim E' +E'' $ or $ L=2E'$ for
some genus 1 pencils $|2E'| , |2E''| $ such that $E'E''= 1$. But $C \sim 3(E_1
+E_2)$ and $gon(C) = M'L'= CL' - L'^2 \geq 4$ .

If $C$ is obtained by Steiner costruction then $C$ contains two base points of
pencil $|L|$ . Since $C$ is very ample there is a smooth curve $C' \in |C|$
which does not contain base points of any pencil $|L'| = |E' + E''| $ ,
because we have at most countable number of such pencils. Now we claim that
$gon(C') \geq 5$. Indeed , by proposition~\ref{PPP} (i) , (iii) we have
$gon(C') = M'L' \geq 4$ and equality occurs  , as we have seen above , only if
$L'^2 =2$ and $C'$ is obtained by Steiner construction using $|L'| , |M'(
\lambda ) |$ but this is not the case , because $C' $ does not contain base
points of linear system $|L'|$.
\newtheorem{G}{1.}
\begin{G}
   A general curve in the linear system $|3(E_1 + E_2 ) |$  has gonality $\geq$
5 and Clifford index $\geq$ 3. There is a linear subsystem $V \subset |C|$ of
codimension 2 such that every smooth curve in $V$ has gonality 4 and Clifford
index 2. Moreover, gon($|C|$) = 4 \\
\end{G}

The subsystem $V \subset |C| $ consists of the curves which contain two base
points of the pencil $|L| = |E_1 + E_2 |$ and therefore is obtained by Steiner
construction  using $|L|$ , $|M|$.

Similarly one can obtain the following:
\newtheorem{GG}{2.}
\begin{GG}
A general curve in the linear system $|3E_1 + 4E_2 ) | $ has gonality 6 and
Clifford index 4 . There is a linear subsystem $V \subset |C|$ of codimension 2
such that every smooth curve in $V$ has gonality 5 and Clifford index 3.
Moreover gon($|C|$) =5.
\end{GG}
In this case $M=3E_1 + 2E_2$ , $L=E_1 + E_2$. And the curve $C \in |M+L|$ which
is obtained by Steiner construction using $|L|$ ,$|M( \lambda) | \subset |M| $
has gonality 5 and Clifford index 3. But $C$ is very ample
and  therefore a general curve of the linear system  $|C|$ does not contain
base points of any pencil $|L'|$ . The proof is similar to the proof in the
example above  and we omit the details.

\paragraph{References} : \\

 [B] F.A. Bogomolov "Holomorphic tensors and vector bundles on projective
varieties" Izvestya of Ac.of Sc. USSR ser. math. v.42 N6 pp.1227-1287. \\

[CD] F. Cossec ; I. Dolgachev "Enriques Surfaces 1", Birkh \"{a}user 1989. \\

[GL] M. Green ; R. Lazarsfeld "Special divisors on curves on K3 surfaces" \ \ \
\ \ \ \ Invent. Math. 89, 357-370, 1989. \\

[ELMS] D. Eisinbud , H. Lange , G. Martens , F.-O. Schreyer "The Clifford
dimension of a projective curve ." , preprint. \\

[P] G. Pareschi "Exceptional linear systems on curves on Del Pezzo surfaces" ,
Math. Ann. 291 , 17-38 (1991) . \\

[T] A.N. Tyurin "Cycles, curves  and vector bundles on algebraic surfaces."
Duke Math.J. 54, 1-26,(1987). \\

Department of geometry and topology,

Faculty of mathematics,

Vilnius university,

Naugarduko g.24,

2009 Vilnius, Lithuania.

\end{document}